\documentclass[11pt,a4paper]{article}
\usepackage{amssymb,latexsym}
\usepackage{a4}
\def\eqnreset{\setcounter{equation}{0}}
\def\eqsection#1{\section{#1}\eqnreset}

\newtheorem{Thm}{Theorem}[section]
\newtheorem{Defi}[Thm]{Definition}
\newtheorem{Cor}[Thm]{Corollary}
\newtheorem{Lemma}[Thm]{Lemma}
\newtheorem{Prop}[Thm]{Proposition}
\newtheorem{Rem}[Thm]{Remark}
\newtheorem{Conj}[Thm]{Conjecture}
\newtheorem{Prelim}[Thm]{Preliminary}

\newenvironment{thm}[0]{\begin{Thm}\noindent}%
{\end{Thm}}
{\end{Defi}}
\newenvironment{cor}[0]{\begin{Cor}\noindent}%
{\end{Cor}}
\newenvironment{lemma}[0]{\begin{Lemma}\noindent}%
{\end{Lemma}}
\newenvironment{prop}[0]{\begin{Prop}\noindent}%
{\end{Prop}}
\newenvironment{rem}[0]{\begin{Rem}\noindent\rm}%
{\end{Rem}}
{\end{Conj}}
{\end{Prelim}}

\def\proof{\par\noindent{\it Proof.}{\ }{\ }}
\def\qed{~\hfill$\square$\medbreak}

\def\naam#1{\label{#1}}
\def\refer#1{\ref{#1}}

\def\bib#1{\cite{#1}}

\def\text#1{\;\;\;\;{\rm \hbox{#1}}\;\;\;\;}
\def\qquad{\quad\quad}

\def\minspace{\vspace{-2mm}}
\def\itema{\item[{\rm (a)}]}
\def\itemb{\item[{\rm (b)}]}
\def\itemc{\item[{\rm (c)}]}

\def\msy#1{{\mathbb #1}}
\def\C{{\msy C}}
\def\N{{\msy N}}
\def\Z{{\msy Z}}
\def\R{{\msy R}}
\def\D{{\msy D}}

\def\ga{\alpha}
\def\gb{\beta}

\def\gf{\varphi}

\def\gl{\lambda}

\def\gs{\sigma}

\def\gD{\Delta}

\def\gS{\Sigma}

\def\fa{{\mathfrak a}}
\def\fb{{\mathfrak b}}
\def\fg{{\mathfrak g}}
\def\fh{{\mathfrak h}}

\def\fk{{\mathfrak k}}

\def\fm{{\mathfrak m}}
\def\fn{{\mathfrak n}}

\def\fp{{\mathfrak p}}
\def\fq{{\mathfrak q}}
\def\fs{{\mathfrak s}}

\def\to{\rightarrow}
\def\Re{\mathrm{Re}\,}

\def\inp#1#2{\langle#1\,,\,#2\rangle}
\def\hinp#1#2{\langle#1\,|\, #2\rangle}
\def\Ad{\mathrm{Ad}}
\def\End{\mathrm{End}}
\def\Hom{\mathrm{Hom}}

\def\ad{\mathrm{ad}}
\def\after{\,{\scriptstyle\circ}\,}

\def\tr{\mathrm{tr}\,}

\def\iq{{\mathrm q}}
\def\iC{{\scriptscriptstyle \C}}

\def\cA{{\mathcal A}}

\def\cC{{\mathcal C}}
\def\cD{{\mathcal D}}

\def\cH{{\mathcal H}}

\def\cP{{\mathcal P}}
\def\cR{{\mathcal R}}

\def\cU{{\mathcal U}}

\def\cW{{\mathcal W}}

\def\col{\,:\,}
\def\faq{\fa_\iq}
\def\faqd{\fa_\iq^*}
\def\faqdc{\fa_{\iq\iC}^*}
\def\fadc{\fa_{\iC}^*}
\def\fad{\fa^*}
\def\nC{C^\circ}
\def\oC{{}^\circ \cC}
\def\oCtau{\oC(\tau)}
\def\nE{E^\circ}
\def\Par{{\cP}}

\def\Parsigma{\Par_\sigma}
\def\Parmin{\Parsigma^{{\rm min}}}

\def\Aq{A_\iq}
\def\supp{\mathop{\rm supp}}
\def\Lie{\mathop{\rm Lie}}

\def\Mfu{\widehat M_{\mathrm{fu}}}
\def\Mps{\widehat M_{\mathrm{ps}}}
\def\ev{\mathrm{ev}}
\def\Ind{\mathrm{Ind}}

\def\Cartan{\theta}
\def\spX{{\mathrm X}}
\def\Ci{C^\infty}
\def\Vtau{V_\tau}
\def\minparabs{\cP_\gs^{\rm min}}
\def\Aqp{A_\iq^+}
\def\spXp{\spX_+}
\def\GL{\mathrm{GL}}
\def\embeds{\hookrightarrow}
\def\too{\longrightarrow}
\def\reg{\mathrm{reg}}
\def\dotvar{\, \cdot\,}
\def\Cminf{C^{-\infty}}

\def\aparabs{{\mathcal P}}
\def\faPdc{\fa_{P\iC}^*}
\def\faPd{\fa_P^*}

\def\faP{\fa_P}
\def\faQ{\fa_Q}
\def\sepbar{\,|\,}
\def\starP{{}^*P}
\def\starQ{{}^*Q}
\def\gL{\Lambda}

\def\starC{{}^*C}
\def\psgpR{S}
\def\cAtwoP{\cA_{2,P}}

\def\replacefbythis{g}
\def\bs{\backslash}
\begin{document}
\title{Polynomial estimates for c-functions\\
on reductive symmetric spaces}
\author{E.~P.~van den Ban and H.~Schlichtkrull}
\date{1/11-2010}
\maketitle
\begin{abstract} The c-functions, related to a
reductive symmetric space $G/H$ and a fixed
representation $\tau$ of a maximal compact
subgroup $K$ of $G$, are shown to satisfy
polynomial bounds in imaginary directions.  \end{abstract}
\section*{Introduction}
Let $X=G/H$ be a reductive symmetric space and
$K\subset G$ a maximal compact subgroup. The harmonic
analysis for $K$-finite functions on $G/H$ has
been developed in the past 3 decades (see for example
\cite{FJdiscr},
\cite{OshMat},
\cite{Delpl},
\cite{DelICM},
\cite{BSpl1},
\cite{BSpl2},
\cite{BSpw}
and further references mentioned there)
as a generalization of
the work of Harish-Chandra \cite{HC3}.
Associated to each finite dimensional representation
$(\tau,V_\tau)$ of $K$, and to each cuspidal
$\sigma$-parabolic subgroup $P$ of $G$,
there is a fundamental family
of functions on $G/H$, called Eisenstein integrals,
on which the spectral decomposition of $L^2(G/H)$ is based.
The asymptotic behavior of these functions,
when the space variable $x\in X$ tends to infinity in various directions,
is described through $c$-functions, analogous to
Harish-Chandra's generalized $c$-functions.
The $c$-functions, which play a profound
role in the harmonic analysis, are  matrix-valued
meromorphi functions of the complex spectral parameter $\lambda$.

It is convenient to normalize the Eisenstein
integrals so that the $c$-function corresponding to
a particular asymptotic direction is 1. The asymptotic
behavior
of the normalized Eisenstein integral is then described
by the normalized $c$-functions, which are
quotients of the unnormalized $c$-functions. They are denoted
$\nC(s\col\gl)$.
In this paper we shall establish a polynomial
bound for these normalized $c$-functions, when
$\gl$ tends to infinity in imaginary
directions, for the case where the parabolic subgroup $P$ is
$\sigma$-minimal.
More precisely we prove that
if $\gl$ varies such that its real part stays bounded,
then there exists a polynomial $q(\gl)$
such that the entries of the product
$q(\gl)\nC(s\col\gl)$ are polynomially bounded in $\gl$.
Obviously, the polynomial $q$ serves
to eliminate singularities. Our result includes the
information about this polynomial that it is a product of
first order polynomials given by roots.

The result which we shall prove was stated without proof as
Theorem 10.1 in \cite{BSpw}. It was then used to obtain a
characterization of the Paley-Wiener space for $G/H$, which
differs from the one established in \cite{BSpw}, and which
had previously been conjectured in \cite{BSmc}.
The essential difference between the two characterizations
is that in \cite{BSmc} the growth conditions are only
required in a single Weyl chamber.
Every
Fourier transform $\varphi$ of a compactly supported function
of type $\tau$ satisfies the transformation property
$\varphi(s\gl)=\nC(s\col\gl)\varphi(\gl)$ for each
element $s$ in the Weyl group. It  follows from
the polynomial estimate for $\nC(s\col\gl)$ that the
growth estimates of Paley-Wiener type hold for all
$\gl$ if they hold on a single chamber.

In the special case of a Riemannian symmetric space
and the trivial $K$-type, the $c$-function is explicitly known
from the formula of Gindikin and Karpelevic,
by which it is expressed as a product of quotients of
gamma functions (see for example \cite{Helgaga}, p.\ 447).
In this case our normalized $c$-functions
are the quotients $c(s\gl)/c(\gl)$, where $s$ belongs to
the Weyl group, and the polynomial estimate follows
easily from known properties of the gamma function.
More generally, when the group $G$ is considered as a symmetric
space for $G\times G$, a similar argument can be carried out for
the quotients which normalize
Harish-Chandra's generalized $c$-functions for this case.
Here one employs Wallach's result \cite{Walcfunc}, that the matrix entries of
the generalized $c$-functions are quotients of
gamma functions.

As in Harish-Chandra's case, the $c$-functions for a $K$-type $\tau$ on
a reductive symmetric space $G/H$ are closely
related to standard intertwining operators for the principal series.
These operators satisfy estimates of polynomial type.
On the subspace of $K$-finite vectors this follows for example from
\cite{KS2}, and in general it follows from
the functional equation of Vogan and Wallach \cite{VWfe}.
In the first part of the paper, Section \ref{s: estimates operators},
we establish an estimate of this type for standard intertwining operators,
with refined information about the polynomial factor that
governs singularities. The estimate is derived by using the functional
equation of \cite{VWfe}.
As we believe this result to be of independent interest,
we have written the exposition of this part without reference
to the theory of reductive symmetric spaces.

The main result in the article, the  estimate
of $c$-functions, is stated
in  Thm.\ \ref{t: ThmSB} in Section~\ref{s: eisenstein integrals}.
In this section we also introduce the
necessary background notation related to
reductive symmetric spaces.
The proof is given in Sections
\refer{s: proof of estimates}-\refer{s: Esteval}.
The main ingredients are the previous estimate for
the standard intertwining operators, a similar estimate for a map which parametrizes $H$-fixed distribution vectors in the minimal
principal series of $G/H,$ and finally a similar estimate for the inverse of this parametrizing map, the
so-called evaluation map.

\eqsection{Estimates for intertwining operators}
\label{s: estimates operators}
Let $G$ be a real reductive group of the Harish-Chandra class,
see \cite{VarHA},
and let $K$ be a maximal compact subgroup of $G,$ with associated
Cartan involution $\Cartan.$ The differential
of this involution is an involution of the Lie algebra $\fg$ and denoted
by the same symbol. As usual, we adopt the convention to denote Lie groups
by Roman capital letters, and the associated Lie algebras
by the corresponding German lower case letters.
Let $\fg = \fk \oplus \fp$ be the decomposition
into a direct sum of the $+1$ and $-1$ eigenspaces
of $\Cartan.$ We extend the Killing form of $[\fg, \fg]$ to an $\Ad(G)$-invariant symmetric bilinear
form $B$ on $\fg,$ for which
$\langle X,Y\rangle:= - B(X, \Cartan Y)$ is a positive definite inner product.

We assume that $\fa_0$ is a maximal abelian subspace
of $\fp.$
Let $\aparabs$ be the set
of parabolic subgroups of $G$ containing $A_0: = \exp \fa_0.$
Then each $P \in \aparabs$ has a Langlands decomposition
of the form $P = M_P A_P N_P,$ with $\fa_P \subset \fa_0.$
If $L$ is a group of the Harish-Chandra class, we denote
by $\cR(L)$ the set of (equivalence classes of) continuous
representations $\xi$ of $L$ in a Hilbert space with the following properties:
\begin{enumerate}
\itema
$\xi$ has an infinitesimal character;
\minspace
\itemb
the space of $C^\infty$-vectors for $L$ equals
the space of $C^\infty$-vectors for a maximal
compact subgroup $K_L$ of $L.$
\end{enumerate}
Note that, in (b) we need only require identity of vector spaces,
since the closed graph theorem then guarantees that the
Fr\'echet spaces of smooth vectors are equal.
As all maximal compact subgroups of $L$ are conjugate,
condition (b) is equivalent to the similar
condition for any fixed maximal compact subgroup of $L.$ We note that
all finite-dimensional irreducible representations of $L$
belong to $\cR(L)$. Likewise, the irreducible
unitary representations belong to $\cR(L)$, see
\bib{Wall2} p. 3.

Let $P \in \aparabs $  and $\xi \in \cR(M_P).$ We denote by $\cH_\xi$ the Hilbert space
in which $\xi$ is continuously realized, and by $\cH_\xi^\infty$ the space
of $C^\infty$ vectors for $M_P.$ Let $\gl \in \faPdc.$
By $C^\infty(P\col \xi \col \gl)$ we denote
the space of $C^\infty$-functions $f: G \to \cH_\xi$ transforming according to the
rule
\begin{equation}\label{e: tr rule pr sr}
f(m a n x) = a^{\gl + \rho_P} \xi(m) f(x),
\end{equation}
for all $x \in G$ and $(m,a,n) \in M_P \times A_P \times N_P.$
Here $\rho_P\in\faPd$ is $1/2$ times
the trace of the adjoint
action of $\fa_P$ on $\fn_P$.
Observe that every map in this space has values in $\cH_\xi^\infty.$
The space $C^\infty(P\col \xi \col \gl)$
is equipped with a Fr\'echet topology in the usual way. The restriction
of the right regular representation to this space defines a continuous
representation of $G,$ which we denote by $\pi_{P, \xi, \gl}.$

We denote by $C^\infty(K \col \xi)$ the space
of $C^\infty$ functions $\gf: K \to \cH_\xi$ transforming according to the rule
\begin{equation}\label{e: tr rule for K}
\gf(mk) = \xi(m) \gf(k),\qquad (k \in K, \, m \in M_P \cap K).
\end{equation}
This space is also equipped with the usual Fr\'echet topology.
It follows from condition (b) in the definition of $\cR(M_P)$
that each function $\gf \in C^\infty(K \col \xi)$ has values in
$\cH_\xi^\infty.$ Since $G=PK$ it is
readily seen that the restriction map $f \mapsto f|_K$
defines a topological linear isomorphism from $C^\infty(P\col \xi \col \gl)$
onto $C^\infty(K \col \xi).$ Accordingly, the representation
$\pi_{P,\xi,\gl}$ may be identified with a continuous representation of
$G$ in $\Ci(K\col\xi),$ which we shall denote by the same symbol.
This is called the compact picture of the principal series
representation.

The results of this section are based on the Vogan-Wallach functional
equation for the standard intertwining operator. This equation has
been established for a class of groups different from the
Harish-Chandra class.
This presents no serious complication, but we have to be somewhat
careful with connected components (see the proof of
Thm. \ref{t: functional equation}).

Let us recall the functional equation in a form
that is  suitable for us. If $P\in \aparabs$
we write $\gS(P, \faP)$ for the set of
$\faP$-weights in $\fn_P$, and we denote by
$\bar P\in\aparabs$ the opposite parabolic subgroup
with $\fa_{\bar P}=\faP$ and $\gS(\bar P,\fa_{\bar P})=-\gS(P,\faP)$.
If $P,Q\in\aparabs$ with $\faP=\faQ$,
then for $R \in \R$ we define the set
\begin{equation}
\label{e: defi faPd Q sepbar P}
\faPd(Q\sepbar P, R): = \{\gl \in \faPdc \mid \inp{\Re \gl}{\ga} > R, \;\;\;\forall \ga
\in \gS(\bar Q, \faP) \cap \gS(P, \faP) \;\}.
\end{equation}
Let now $\xi \in \cR(M_P).$ Then there exists a constant $r_\xi > 0$
such that for
$f \in \Ci(P\col \xi\col \gl)$ and $\gl \in \faPd(Q \sepbar P, r_\xi)$
the integral
\begin{equation}
\naam{e: defining integral for A}
A(Q\col P \col \xi\col \gl)f (x) = \int_{N_Q \cap \bar N_{P}}
f(n x) \; dn, \qquad (x \in G)
\end{equation}
is absolutely convergent and defines a function in $C^\infty(Q\col \xi \col
\gl).$  Here $dn$ is a suitably normalized Haar measure on
the nilpotent group $N_Q \cap \bar N_{P}$ whose precise
normalization
is immaterial for our discussion.

For $\gl \in \faPd(Q \sepbar P, r_\xi),$ the operator
\begin{equation}
\naam{e: standard intertwiner}
A(Q\col P \col \xi \col \gl):\;\;
C^\infty(P\col \xi\col \gl) \to C^\infty(Q \col
\xi \col \gl)
\end{equation}
is continuous linear. It intertwines the representations
$\pi_{P, \xi, \gl}$ and $\pi_{Q, \xi, \gl}$ and is called
the standard intertwining operator between these representations. Finally, for all
$f \in \Ci(K\col \xi)$
the map $\gl \mapsto [A(Q \col P \col \xi \col \gl) f]|_K$
is a holomorphic $C (K\col \xi)$-valued map,
and there exists a constant $C > 0 $ such
that the following estimate is valid, for all $\gl \in \faPd(Q\sepbar P, r_\xi)$ and
all $f \in \Ci(P\col \xi \col \gl):$
\begin{equation}
\naam{e: estimate sup norm of A f}
\|A(Q \col P \col \xi \col \gl) f\|_0 \leq C \|f \|_0.
\end{equation}
Here $\|\dotvar\|_0$ denotes the supremum norm over $K$ of a
function $G \to \cH_\xi.$
All these facts are  found
in \bib{VWfe} or \bib{Wall2},
and they are easily seen to be
valid also for all groups of the Harish-Chandra class.
If we combine the above
estimate with the equivariance of $A(Q \col P \col \xi \col \gl)$
we readily infer that
$\gl \mapsto A(Q \col P \col \xi \col \gl) f$ is  a holomorphic
$\Ci(K\col \xi)$-valued
function on $\faPd(Q|P, r_\xi),$ for every $f \in \Ci(K\col \xi).$

\begin{thm} {\rm (Vogan--Wallach \bib{VWfe})\ }
\naam{t: functional equation}
Let $P \in \aparabs$ and $\xi \in \cR(M_P).$
There exist non-trivial polynomial functions $b_\xi: \faPdc \to \C$
and $D_\xi: \faPdc \to U(\fg)^{K}$ such that for every
$f \in C^\infty(P\col \xi\col \gl)$ and all $\gl \in \faPd(\bar P \sepbar P , r_\xi)$
the following equality  is valid:
\begin{eqnarray}
\lefteqn{b_\xi(\gl)\,  A(\bar P \col P \col \xi \col \gl)\, f }\nonumber\\
\naam{e: functional equation}
&=&
 A(\bar P\col P \col \xi\col \gl + 4 \rho_P)
\,\pi_{P, \xi, \gl + 4 \rho_P}(D_\xi(\gl))\,f.
\end{eqnarray}
\end{thm}

\proof
If $G$ is connected, then  it is readily seen that $G$ belongs
to the class of groups considered in \bib{VWfe}; see \bib{BWcc},
Sect.\ 0.3.1,
for its definition.
The result then follows from \bib{VWfe}, Thm.\ 1.5
(see also \bib{Wall2}, Thm.\ 10.1.5).

Now assume that
$G$ is a general group of the Harish-Chandra class.
Then the identity component
$G_e$ is of the Harish-Chandra class. Moreover, the representation
 $\xi_e:= \xi|_{M_P \cap G_e}$ has the same infinitesimal character as $\xi.$
The smooth vectors for $\xi_e$  are the same as those for $\xi,$
hence as those for $\xi|_{M_P \cap K},$ which in turn are the smooth vectors for
the restriction of $\xi_e$ to the open subgroup $M_P \cap K_e$ of $M_P \cap K.$
It follows that $\xi_e \in \cR(M_P \cap G_e).$ Therefore, the
above result is valid for the
data $G_e,$ $P \cap G_e,$ $M_P \cap G_e$ and $ \xi_e.$
It follows that there exist polynomial maps
 $b_\xi: \faPdc \to \C$ and  $D_\xi: \faPdc \to U(\fg)^{K_e}$ with
\begin{eqnarray*}
\lefteqn{
b_\xi(\gl)A(\bar P \cap G_e
\col P \cap G_e  \col \xi_e \col \gl) (f|_{G_e})}\\
&=&
A(\bar P \cap G_e \col P \cap G_e \col \xi_e \col \gl + 4 \rho_P)
\pi_{P \col G_e, \xi_e,  \gl + 4 \rho_P}(D_\xi(\gl))(f|_{G_e})
\end{eqnarray*}
for all $f \in \Ci(P\col \xi \col \gl)$ and all
$\gl$ in the set $\faPd(\bar P \sepbar P , r_\xi)$,
where the intertwining operator on the left is given by
an absolutely convergent integral. In fact, since the integrals
for $G$ and $G_e$ are identical, we have
$$A(\bar P \cap G_e \col P \cap G_e  \col \xi_e \col \gl) (f|_{G_e})
=[A(\bar P\col P\col\xi\col\gl)f]|_{G_e}
$$
for all $\faPd(\bar P \sepbar P , r_\xi)$. Similarly
\begin{eqnarray*}
\lefteqn{
A(\bar P \cap G_e \col P \cap G_e \col \xi_e \col \gl + 4 \rho_P)
\pi_{P \col G_e, \xi_e,  \gl + 4 \rho_P}(D_\xi(\gl))(f|_{G_e})}\\
&\qquad \qquad \qquad =&
[A(\bar P\col P \col \xi\col \gl + 4 \rho_P)
\,\pi_{P, \xi, \gl + 4 \rho_P}(D_\xi(\gl))\,f]|_{G_e}.
\end{eqnarray*}
It thus follows that
\begin{eqnarray}
\lefteqn{b_\xi(\gl)\,  [A(\bar P \col P \col \xi \col \gl)\, f ](g)}\nonumber\\
&=&
 [A(\bar P\col P \col \xi\col \gl + 4 \rho_P)
\,\pi_{P, \xi, \gl + 4 \rho_P}(D_\xi(\gl))\,f] (g).
\label{e: functional equation at g}
\end{eqnarray}
for all $g\in G_e$.
Since $G=M_{P}G_e$ the identity (\ref{e: functional equation at g})
then follows for $g\in G$
by means of (\ref{e: tr rule pr sr}).

It only remains to be seen that (\ref{e: functional equation})
can be arranged with $D_\xi: \faPdc \to U(\fg)^{K}$
instead of $D_\xi: \faPdc \to U(\fg)^{K_e}$. Let $k\in K$.
By application
of  (\ref{e: functional equation at g}) with $R_{k^{-1}}f$ and $gk$
in place of $f$ and $g$, respectively, we derive easily
from the intertwining property of $A$
that
\begin{eqnarray*}
\lefteqn{b_\xi(\gl)\,  A(\bar P \col P \col \xi \col \gl)\, f }\\
&=&
 A(\bar P\col P \col \xi\col \gl + 4 \rho_P)
\,\pi_{P, \xi, \gl + 4 \rho_P}(\Ad(k)D_\xi(\gl))\,f.
\end{eqnarray*}
The desired identity (\ref{e: functional equation})
now follows with $D_\xi(\gl)$ replaced by the integral over
$K$ of $ \Ad(k)D_\xi(\gl)$. It is easily seen that this
integral is again a polynomial in $\gl$, and it is non-trivial
because the left side of  (\ref{e: functional equation}) is
non-trivial for some $f$.
\qed

The functional equation leads to meromorphic continuation
of the standard intertwining operator with suitable estimates.
We will cast these estimates
in a form which is suitable for the rest of this paper.
We need the following preparations.

For $s\in\N=\{0,1,\dots\}$
we denote by $C^s(K \col \xi)$ the space
of $C^s$ functions $\gf: K \to \cH_\xi$ transforming
according to the rule
(\ref{e: tr rule for K}). As usual we write $C$ for $C^0$.
As $K$ is compact, it follows that $C^s(K\col\xi)$
can be equipped with a Banach norm
$\|\,\cdot\,\|_s$ as follows. Let $X_1,\dots,X_m$ be a basis
for $\fk$, and let $X^n=X_1^{n_1}\cdots X_m^{n_m}\in\cU(\fk)$
for a multi-index $n\in\N^m$, then
\begin{equation}\label{d: s-norm}
\|f\|_s:=\sum_{|n|\leq s} \|R_{X^n}f\|_0.
\end{equation}
It is easily seen that the right regular action of $K$ on $C^s(K\col\xi)$
is bounded with respect to this norm.

For every $s \in \N,$ it follows from
(\ref{e: estimate sup norm of A f}) and the intertwining property
of $A(Q\col P\col\xi\col\gl)$ that there exists a constant
$C>0$ such that
the following estimate is valid, for all $f \in \Ci(K\col \xi)$
and $\gl \in \faPd(Q\sepbar P, r_\xi)$
\begin{equation}
\naam{e: estimate r norm of A f}
\|A(Q \col P \col \xi \col \gl) f \|_s \leq C \|f\|_s.
\end{equation}
It follows that the intertwining operator extends to a bounded linear
endomorphism of the Banach space $C^s(K\col \xi),$ for $\gl$ in the
indicated region. Moreover,
the operator norm of this extension is uniformly bounded in $\gl.$
By an easy application of the Cauchy integral formula, it
follows that the
intertwining operator is holomorphic in the variable
$\gl \in \faPd(Q|P, r_\xi),$
as a function with values in the Banach space $B(C^s(K \col \xi))$ of bounded linear endomorphisms, equipped
with the operator norm.

\begin{lemma}
\naam{l: estimate pi D}
Let $D \in U(\fg)$ be an element of order at most $d$.
Then there exists a constant $t \in \N,$ and for every $s \in \N$
a constant $C >0,$ such that
\begin{equation}
\naam{e: estimate pi D}
\| \pi_{P, \xi, \gl}(D)f \|_s \leq C \,  (1 + |\gl|
)^d\, \|f\|_{s + t},
\end{equation}
for all $f \in \Ci(K\col \xi)$ and all $\gl\in\faPdc$.
\end{lemma}

\proof
It clearly suffices to prove the result for $d=1$ and
$D\in\fg$. We use the decomposition
$\fg = \fn_P \oplus \fa_P\oplus (\fm_{P}\cap\fp) \oplus \fk$.
Since $\fk$ is $\Ad(K)$-invariant, it follows that
$$\fg = \Ad(k^{-1})\Big(\fn_P \oplus \fa_P\oplus (\fm_{P}\cap\fp)
\Big)\oplus \fk$$
for each $k\in K$.
Accordingly, we write
$$D= \Ad(k^{-1})(U_k + V_k +W_k) + Z_k. $$
The elements $U_k\in\fn_P$,
$V_k\in\fa_P$, $W_k\in\fm_P\cap\fp$ and $Z_k\in\fk$ all
depend smoothly on $k\in K$.

If $f \in \Ci(K\col \xi)$ and $\gl \in \faPdc,$ then by $f_\gl$ we denote
the unique function in $\Ci(P\col \xi \col \gl)$ that restricts to $f$ on $K.$
Then
\begin{eqnarray}
\lefteqn{
[\pi_{P, \xi, \gl}(D)f](k)}\nonumber\\
&=&
[L_{-U_k-V_k-W_k}f_\gl](k)+
[R_{Z_k}f](k)\nonumber\\
&=& (\gl+\rho_P)(V_k)f(k)+\xi(W_k)(f(k))+[R_{Z_k}f](k).
\label{three terms}
\end{eqnarray}
In order to estimate $\|\pi_{P, \xi, \gl}(D)f\|_s$
we use the definition (\ref{d: s-norm}).
Let $n$ be a multi-index with $|n|\leq s$, and
apply $X^n\in\cU(\fk)$ to (\ref{three terms}),
as a function of $k$. We consider the three terms separately.

Since $k\mapsto V_k$ is smooth, the derivatives of
$k\mapsto (\gl+\rho_P)(V_k)$ up to order $s$
are bounded uniformly by a constant times $1+|\gl|$.
By the Leibniz rule, it then follows that
$$\|\,R_{X^n}[k\mapsto (\gl+\rho_P)(V_k)f(k)]\,\|_0\leq
C(1+|\gl|)\|f\|_s$$
for some constant $C>0$.

Let $W \in \fm_{P}.$ Then $\xi(W)$ maps
$\cH_\xi^\infty$ continuously
into itself. By assumption (b) in the definition of $\cR(M_P)$
there exists a finite subset $F \subset U(\fm_P \cap \fk)$
such that
$$
\|\xi(W) a \|_{\xi} \leq \max_{u\in F} \|\xi(u) a\|_{\xi}
, \qquad (a \in \cH_\xi^\infty).
$$
This in turn implies that there exist constants
$r\in \N$ and $C > 0$
such that, for all $f \in \Ci(K\col \xi)$ and $k\in K$,
\begin{eqnarray*}
\|\xi(W) (f(k))\|_{\xi}
&  \leq & \max_{u \in F} \| \xi(u) (f(k))\|_{\xi}  \\
&=&  \max_{u \in F} \| [L_{u^\vee}f](k))\|_{\xi}  \\
&=&  \max_{u\in F} \|[ R_{\Ad(k)^{-1} u}f ](k) \|_{\xi} \\
&  \leq  & C \|f\|_{r}.
\end{eqnarray*}
Since the right action of $K$ commutes with the application of
$\xi(W)$ in this expression, we conclude that
the derivatives of $k\mapsto \xi(W)(f(k))$ up to
order $s$ are bounded by $C\|f\|_{r+s}.$
It follows that
$$\|\,R_{X^n}[k\mapsto \xi(W_k)(f(k))]\,\|_0\leq
C\|f\|_{r+s},$$
for some constant $C>0$.

Finally, it is clear that
$$\|\,R_{X^n}[k\mapsto R_{Z_k}f(k)]\,\|_0\leq
C\|f\|_{s+1}$$
for some constant $C>0$.
The lemma now follows with $t=\max\{r,1\}$.
\qed

We proceed by proving an estimate from which we shall derive
a more general estimate in Theorem \refer{t: estimate A on smooth func},
by product decomposition.
We need the following simple observation.

\begin{lemma}\naam{positive inner product}
$\inp{\rho_P}{\ga}>0$ for all $\ga \in \gS(P, \faP).$
\end{lemma}

\proof
Let $P_0 \subset P$ be a minimal parabolic
subgroup containing $A_0.$
Then it is well known that
$\inp{\rho_{P_0}}{\ga} > 0$ for all $\ga \in \gS(P_0, \fa_0).$
The latter set is a positive system for $\gS(\fg, \fa_0)$; let $\gD(P_0)$ be
the subset of simple roots.

Put $\fa_{M_P}:= \fa_0 \cap \fm_P$ and let $\rho_{M_P} \in \fa_{M_P}^*$ be the
$\rho$ of the minimal parabolic subgroup $P_0 \cap M_P$ of $M_P.$
Then $\rho_{P_0} = \rho_{M_P} + \rho_P,$ by \cite{KS2}, eq.~(1.11),
and this decomposition is compatible with the decomposition
$\fa_0 = \fa_{M_P} \oplus \fa_P.$

Let now $\ga \in \gD(P_0).$ If $\ga$ restricts to zero on $\fa_P,$ then
$\inp{\rho_P}{\ga} = 0.$ On the other hand, if $\ga$ does not restrict to zero
on $\faP,$ then $\inp{\rho_{M_P}}{\ga} \leq 0$ and we see that
$\inp{\rho_P}{\ga} > 0$ for such $\ga.$ From these properties of the simple
roots it follows that
for each $\ga \in \gS(P_0, \fa_0)$ the inner product $\inp{\rho_P}{\ga}$
is positive as soon as $\ga|_{\fa_P} \neq 0.$ The non-zero restrictions
$\fa|_{\fa_P}$ constitute $\gS(P, \faP),$ and since $\rho_P= 0$
on $\fa_{M_P},$ the claim follows.
\qed

\begin{cor}
\naam{c: mero cont of A with est}
Let $P \in \aparabs$ and $\xi \in \cR(M_P).$ Then for every $f \in \Ci(K\col \xi)$
the $\Ci(K\col \xi)$-valued function $\gl \mapsto A(\bar P \col P \col \xi \col \gl)f$
extends meromorphically to $\faPdc.$

For every $R \in \R$ there exists a polynomial function $q: \faPdc \to \C$ and  constants $N \in \N$ and $r \in \N$
with the following properties.
\begin{enumerate}
\itema
For every $f \in \Ci(K\col \xi)$ the $\Ci(K\col \xi)$-valued meromorphic function
$\gl \mapsto q(\gl) A(\bar P \col P \col \xi \col \gl) f$
is regular on $\faPd(\bar P | P, R).$
\itemb
For every $s \in \N$ there exists a constant $C >0$
such that, for all $f \in \Ci(K\col \xi)$
and all $\gl\in\faPd(\bar P | P, R),$
\begin{equation}
\naam{e: estimate q A first}
\|q(\gl) A(\bar P \col P \col \xi \col \gl) f\|_s \leq C (1 + |\gl|)^N \|f\|_{s + r}.
\end{equation}
\end{enumerate}
\end{cor}

\proof
We will show that $A(\gl):= A(\bar P \col P \col \xi \col \gl)$
extends meromorphically to $\faPd(\bar P | P , R),$ with the required
estimates, by downward induction on $R\in\R.$
If $R \geq r_\xi$ then $A(\gl)$ depends holomorphically on
$\gl \in \faP(\bar P| P, R),$
and (a) and (b) are valid by (\refer{e: estimate r norm of A f}),
with $q = 1,$ $N = 1$ and $r = 0.$

Let $c$ be the minimum value of
$4\inp{\ga}{\rho_P}$ as $\ga \in \gS(P, \faP).$
It follows from Lemma \ref{positive inner product}
that $c > 0.$
Let $R_1\in\R$. Assuming that
the claimed result has been established for $R = R_1,$
we will show that it is valid for $R = R_1 - c.$
Note that with this choice of $R$ the set
$\faPd(\bar P|P, R) +4 \rho_P $ is contained in
$ \faPd(\bar P| P, R_1).$

Let $q_1:\faPdc \to \C$ be a polynomial and
$N_1,r_1\in\N$ constants, such that
(a) and (b) hold for $R=R_1$,
and let
$b_\xi$ and $D$ be as in Theorem \refer{t: functional equation},
also for $R=R_1$.
Then we will verify (a) and (b) for $R=R_1-c$
with  $q(\gl) := b_\xi(\gl) q_1(\gl + 4 \rho_P).$

We recall that
$D$ is a polynomial on $\faPdc,$ with values in $U(\fg).$ This means that
we can write $D(\gl)$ as a linear combination of finitely many fixed
elements from $U(\fg),$ say of order $\leq d$,
with coefficients that depend polynomially
on $\gl$, say of degree $\leq m$.

Let $f \in \Ci(K\col \xi).$ Then by analytic
continuation, (\refer{e: functional equation})
holds for  $\gl \in \faPd(\bar P|P, R_1).$
It follows that we have the following identity of meromorphic
$C^\infty(K \col \xi)$-valued functions in the variable $\gl \in
\faPd(\bar P|P, R_1):$
\begin{eqnarray*}
q(\gl) A(\gl) f & = &  b_\xi(\gl) q_1(\gl + 4 \rho_P) A( \gl) f \\
&=& q_1(\gl + 4 \rho_P) A(\gl + 4 \rho_P) \pi_{P, \xi, \gl + 4 \rho_P}(D(\gl)) f.
\end{eqnarray*}
The last expression extends holomorphically to
$\gl+4\rho_P\in\faPd(\bar P,|P, R_1),$
hence in particular to $\gl\in\faPd(\bar P | P, R),$
and therefore so does $q(\gl) A(\gl) f.$
The required estimate follows by application
of the induction hypothesis and Lemma \refer{l: estimate pi D}:
\begin{eqnarray*}
\|q(\gl) A(\gl) f\|_s & = &  \| q_1(\gl + 4 \rho_P) A(\gl + 4 \rho_P)
\pi_{P, \xi, \gl + 4 \rho_P}(D(\gl)) f\|_s\\
&\leq& C_1(1+|\gl|)^{N_1}\|
\pi_{P, \xi, \gl + 4 \rho_P}(D(\gl)) f\|_{s+r_1}\\
&\leq& C_2(1+|\gl|)^{N_1+m+d} \|f\|_{s+r_1+t}.
\end{eqnarray*}
\qed

The above result can be refined as follows.
If $V$ is a finite
dimensional real vector
space equipped with a positive definite inner product, and $S\subset V\setminus\{0\}$, we denote by
$\Pi_S(V)$ the set of polynomial functions  on $V_\iC$, which are products
of first order polynomials of the form $\gl \mapsto \inp{\gl}{\ga} -c$ with
$\ga \in S$ and $c\in\C$.
We define $\Pi_{S,\R}(V)$ similarly, but with $c\in\R$.

\begin{thm}
\naam{t: estimate A on smooth func}
Let $P, Q \in \aparabs$ be such that $\fa_P =\fa_Q.$
Let $\xi \in \cR(M_P)$ and $R\in \R.$

There exists a polynomial function
$q \in \Pi_{\gS(\bar Q, \faQ) \cap \gS(P, \faP)}(\faPd)$
such that for every $f \in \Ci(K\col \xi),$ the function
$$
\gl \mapsto q(\gl)A(Q \col P \col \xi\col \gl) f
$$
extends to a holomorphic function
on $\faPd(Q | P, R)$
with values in $C^\infty(K\col \xi),$
for which the following estimates hold:

There exist $r \in \N,$ $N \in \N$
and for every $s \in \N$ a constant $C >0,$ such that
\begin{equation}
\naam{e: estimate q A general}
\|q(\gl) A(Q \col P \col \xi \col \gl) f\|_s \leq C ( 1 + |\gl|)^N \|f\|_{s + r},
\end{equation}
for all  $f \in \Ci(K\col \xi)$
and all $\gl\in\faPd(Q|P,R)$.
\end{thm}

\proof
We first recall
that the operator $A(Q \col P \col \xi \col \gl)$ may be expressed
as a composition of similar operators with
$Q$ and $P$ adjacent,
i.e., the chambers $\fa_P^+$ and $\fa_Q^+$ are adjacent in
$\faP$ or, equivalently, the roots in
$\gS(\bar Q, \faP) \cap \gS(P, \faP)$ are proportional;
see \bib{HC3}, Sect.\ I.2,
\bib{KS2}, Sect.\ 7, and \bib{VWfe} for details.
The result is now readily reduced to the case that $Q$ and $P$ are adjacent.
From now on, we assume this to be the case.

There exists precisely one reduced root $\ga \in \gS(P, \faP)$
such that
$\bar \fn_Q \cap \fn_P = \fn_\ga := \sum_{c>0} \fg_{c\ga}.$
Fix an element $X \in \ker\ga \cap \bar\fa_P^+$
such that $\ga$ is the unique reduced root in
$\gS(P, \faP)$ vanishing on $X.$ Let $\gL\subset \R$ be the collection
of nonnegative weights of $\ad X$ in $\fg.$ The sum
$\fs:= \sum_{\gl \in \gL} \fg_\gl$
of the associated weight spaces is a parabolic subalgebra of $\fg.$
Let $\psgpR  = M_\psgpR  A_\psgpR  N_\psgpR $
be the associated parabolic subgroup in $G$
with the indicated Langlands decomposition. Then
$$
\fm_\psgpR  = \fm_P \oplus \fa_\ga \oplus \fn_\ga \oplus \bar\fn_\ga,
$$
where $\fa_\ga$ is the (one dimensional) orthocomplement
of $\ker \ga$ in $\faP.$  We now observe that
$\starP := P \cap M_\psgpR $ and $\starQ:= Q \cap M_\psgpR $
are (maximal) parabolic subgroups of $M_\psgpR $ with Langlands
decompositions $\starP = M_P A_\ga N_\ga$ and
$\starQ = M_P A_\ga\bar N_\ga.$
In particular, $\starP$ and $\starQ$ are opposite.
From the integral formulas for the standard intertwining operators,
it follows that, for all $f \in \Ci(K\col \xi),$  $k \in K,$
$$
[A(Q \col P \col \xi \col \gl)f](k) =
A(\starQ \col \starP \col \xi \col \gl|_{\fa_\ga})[R_k f|_{M_\psgpR }](e),
$$
for $\gl \in \faPdc$ with $\inp{\Re \gl}{\gb} > c_\xi$ for all
$\gb \in \gS(P, \faP).$ From this it follows by application
of Corollary \refer{c: mero cont of A with est} that the expression
on the left-hand side has a meromorphic extension
as a $\cH_\xi$-valued function of $\gl \in \faPdc.$
Moreover, this function depends on  $\gl$ through the one-dimensional
restriction $\gl|_{\fa_{\ga}}.$ As
$R_k f|_{K_\psgpR }\in \Ci(K_\psgpR \col \xi)$
depends continuously on $k,$ it follows that
$A(Q \col P \col \xi \col \gl)f$ is a meromorphic function
of $\gl|_{\fa_\ga}$ with values in $C(K\col \xi).$

Let now $R \in \R$ and let ${}^* q,$ $N$ and $r$ be associated with
the data $M_\psgpR , \starP, \xi, R$ as in
Corollary \refer{c: mero cont of A with est}.
Put $q(\gl):=
{}^* q (\gl|_{\fa_\ga}).$ Then, clearly, $q \in \Pi_{\{\ga\}}(\faPd).$
Moreover, if $\gl \in \faPd(Q|P, R),$ then
$\gl|_{\fa_\ga} \in \fa_\ga^*(\starQ|\starP , R).$ It follows
that the $C(K\col \xi)$-valued meromorphic function
$\gl \mapsto q(\gl)A(Q\col P\col \xi\col \gl)f$
is regular on $\faPd(Q | P , R)$ and satisfies the estimate
\begin{eqnarray}
\|q(\gl)A(Q\col P\col \xi\col \gl)f\|_0 &\leq &
\sup_{k\in K}{\starC}
(1 + |(\gl|_{\fa_\ga})|\,\|)^N
\|R_k f\|_{C^r(K_\psgpR \col \xi)}
\nonumber \\
\naam{e: estimate sup norm q A}
&\leq & C (1 + |\gl|)^N \|f\|_r,
\end{eqnarray}
with $C > 0$ a constant independent of $f.$
By analytic continuation it follows that the operator
$A(Q\col P\col \xi\col \gl)$ intertwines the representations
$\pi_{P, \xi, \gl}$ and $\pi_{Q, \xi, \gl}$ for generic $\gl \in \faPdc.$
In particular, it is $K$-intertwining, and the estimate
(\refer{e: estimate q A general}) readily follows from
(\refer{e: estimate sup norm q A}). The proof
is now easily completed.
\qed

\begin{rem}\label{r: PiR}
Although we do not need it in the sequel, we note
that if $\xi$ has real infinitesimal character then
Thm.~\ref{t: estimate A on smooth func}
can be improved so that
$q \in \Pi_{\gS(\bar Q, \faQ) \cap \gS(P, \faP),\R}(\faPd)$.
This can be seen from the stated version of the theorem
together with \cite{KS2}, Theorem 6.6,
by means of the argument in
Remark \ref{r: some implies all} below.
\end{rem}

\eqsection{Normalized Eisenstein integrals and $c$-functions}
\label{s: eisenstein integrals}
In the rest of this paper we keep assuming that $G$ is a group of the
Harish-Chandra class unless specified otherwise.
Let $\gs$ be an involution of $G$ and $H$ an open subgroup
of the group $G^\gs$ of fixed points for $\gs.$  Then
$\spX:= G/H$ is a reductive symmetric space of the Harish-Chandra
class.

There exists a Cartan involution $\Cartan$ which commutes
with $\gs.$ The associated maximal compact subgroup is denoted
by $K.$  We have
the decompositions $\fg = \fk \oplus \fp = \fh \oplus \fq$
into direct sums of the $+1$ and $-1$ eigenspaces
of the infinitesimal involutions $\Cartan$ and $\gs,$ respectively. We fix an inner product $\inp{\dotvar}{\dotvar}$
on $\fg$ as in the beginning of Section \ref{s: estimates operators}, subject to the additional requirement
that the involution $\gs$ be symmetric with respect to it.

We fix a maximal abelian subspace $\faq$ of $\fp \cap \fq$ and extend
it to a maximal abelian subspace $\fa_0$ of $\fp.$
Of course the results of the previous
section are available for the present choices of $K$ and $\fa_0.$

We fix a finite dimensional unitary representation $(\tau, \Vtau)$ of
$K$ and consider the following space of smooth $\tau$-spherical
functions on $\spX$:
\begin{equation}
\naam{e: defi spher func}
\Ci(\spX \col \tau): = \{f \in \Ci(\spX,V)\mid
f(kx) = \tau(k)f(x) \;\;\forall x \in
\spX, k \in K\}.
\end{equation}
Let  $\Aq:= \exp \faq.$ Then $G = K \Aq H,$ hence a function
from the space (\refer{e: defi spher func}) is completely
determined by its restriction to $\Aq.$

The restricted root system $\gS$ of
$\faq$ in $\fg$ is a possibly non-reduced root system. The associated collection
of regular points in $\Aq$ is denoted by $\Aq^\reg.$ We recall that the
set $\spXp:= K \Aq^\reg H$ is open dense in $\spX.$

Let $\minparabs$ denote the set of minimal $\gs\Cartan$-stable
parabolic subgroups of $G$ containing $\Aq.$
Then $\minparabs\subset\aparabs$.
{}From \cite{Bps1}, \S 2, we recall that each
$P \in \minparabs$ has a Langlands decomposition of the form
\begin{equation}
\naam{e: Langlands deco P min}
P = MA N_P,
\end{equation}
with $M=M_P$ and $A=A_P$ independent of $P$ and defined as follows.
Let $M_1$ denote the centralizer of
$\faq$ in $\fg,$ then $M_1=MA$ where $M ={}^\circ M_1$
(see \cite{VarHA}, Part II, p.\ 18 for this notation)
and $A = \exp(\mathrm{center}\,(\fm_1) \cap \fp).$
We also recall that
$
\faq\subset \fa  \subset \fa_0.
$
Moreover, the space $\fa$ is $\sigma$-invariant, and
$\fa\cap \fq = \faq,$ so that
\begin{equation}
\naam{e: deco fa}
\fa=(\fa \cap\fh)\oplus \faq.
\end{equation}

For each $P \in \minparabs,$
we denote by $\gS(P)$ the collection of $\faq$-weights
in $\fn_P.$  Then $\gS(P)$ is a positive
system for $\gS;$ the associated open Weyl chamber in
$\Aq$ is denoted by $\Aqp(P).$
By \cite{Bps1}, Lemma 2.8, the map $Q \mapsto \gS(Q)$ is
a bijection from $\minparabs$ onto the collection of positive
systems for~$\gS.$

The image
of the natural embedding $N_K(\faq) / Z_K(\faq) \embeds \GL(\faq)$
equals the Weyl group $W$ of $\gS.$ The image of $N_{K\cap H}(\faq)$
in $W$ is denoted by $W_{K\cap H}.$ We fix a complete
set of representatives $\cW$ of $W/ W_{K\cap H}$ in $N_K(\faq).$
Then for each $Q \in \minparabs$ we have
the disjoint union
$$
\spXp= \cup_{v \in \cW}\;\; K \Aqp(Q) v  H.
$$
Consequently, any function from
(\refer{e: defi spher func}) is completely
determined by its restrictions to $\Aqp(Q)v,$ for $v \in \cW.$

We denote by $\Mfu$ the collection of (equivalence classes)
of finite dimensional irreducible  unitary representations $\xi$ of $M$,
and by $\Mps$ the subcollection of those $\xi$ for which
there exists
$v\in\cW$ such that $\xi$ has a non-trivial
$M\cap vHv^{-1}$-fixed vector.
Of course these representations are trivial on the non-compact factors of $M.$
If $P\in \minparabs,$
then the associated series of induced
representations $\Ind_P^G(\xi \otimes \gl \otimes 1),$ for
$\xi \in \Mps$ and $\gl \in \faqdc$,
is called the minimal principal series
for $\spX,$ see \cite{Bps1},
\S{} 3.

Associated with $\tau$ and this series of representations, the normalized
Eisenstein integral $\nE(P\col \gl)$ is defined as in \cite{BSft}, \S 5.
Let
\begin{equation}\label{defi oCtau}
\oCtau := \oplus_{v \in \cW} \;\;C^\infty(M/M \cap v H v^{-1}\col \tau_M),
\end{equation}
with $\tau_M := \tau|_{M\cap K}.$ This space was  denoted $\cAtwoP$ in
\cite{BSpl1}.
Equipped with the direct sum of the $L^2$-inner products on the summands, it is
a finite dimensional Hilbert space.

The Eisenstein integral is a smooth function on $\spX,$  with values in the space
$\Hom(\oCtau, \Vtau),$ and
depending meromorphically on the parameter $\gl \in \faqdc.$
For each $\psi \in \oCtau,$ and regular $\gl \in \faqdc,$ the function $\nE(P\col \gl)\psi$
belongs to $C^\infty(\spX\col\tau)$ and is
finite for the algebra $\D(\spX)$
of invariant differential operators on $\spX$.
Accordingly, it has a convergent
series expansion on each region $K\Aqp(Q)vH,$ for $Q \in \minparabs$ and
$v \in \cW.$ The normalized $c$-functions
$\nC_{Q|P}(s\col \gl),$ for $s \in W,$ are determined by the leading coefficients in these expansions.
More precisely, each function $\nC_{Q|P}(s \col \dotvar)$ is a
meromorphic function on $\faqdc$
with values in $\End(\oCtau).$  On $i \faqd$ it is regular and unitary,
and for $\gl \in i \faqd$
the top order asymptotic behavior of the Eisenstein integral along $K\Aqp(Q)vH$ is described
by the formula
$$
\nE(P\col \gl \col mav)\,\psi \sim \; \sum_{s \in W} a^{s \gl - \rho_Q}\; [\nC_{Q|P}(s \col \gl)\psi]_v (m),
$$
as $m \in  M,$ and $a \to \infty$ in $\Aqp(Q).$ Here $v\in\cW$ and
the subscript indicates that the $v$-component of an element in the space
(\ref{defi oCtau}) has been taken.

Let $P\in\aparabs$ be fixed, and write
$\nC(s\col \gl): = \nC_{P|P}(s \col \gl).$ The following
result is stated without proof in \cite{BSpw}, Thm.~10.1.
We will give the proof in the  remainder of this article.
The set of polynomials $\Pi_{\gS,\R}(\faqd)$ is defined
above Thm. \ref{t: estimate A on smooth func}.

\begin{thm}
\naam{t: ThmSB}
Let $s\in W$ and let $\omega \subset \faqd$ be compact.
There exists a polynomial $q\in\Pi_{\gS,\R}(\faqd)$ and a number
$N\in\N$ such that
$\gl\mapsto(1+|\gl|)^{-N} q(\gl)\nC(s\col\gl)$ is bounded on $\omega+i\faqd$.
\end{thm}

\eqsection{Polynomial estimates for the $c$-functions}
\naam{s: proof of estimates}
In the proof of Theorem \refer{t: ThmSB}, it will be convenient to not only consider the
functions $\nC(s\col \gl),$ but more generally all
functions $\nC_{Q|P}(s \col \gl),$ for $P,Q \in \minparabs$ and $s \in W$.
In fact, we will establish the following result.

\begin{thm}
\naam{t: main thm}
Let $P,Q\in \minparabs,$ $s\in W$ and let $\omega\subset\faqd$ be compact.
There exists a polynomial function $q\in\Pi_{\gS,\R}(\faqd)$,
a number $N\in\N$ and a constant $C>0$
such that
\begin{equation}
\naam{e: Eqfa}
\| q(\gl)\nC_{Q|P}(s\col\gl) \| \leq C\, (1+|\gl|)^{N},
\end{equation}
for all $\gl \in \omega+i\faqd.$
\end{thm}

\begin{rem}\label{r: some implies all}
Notice that if (\ref{e: Eqfa}) holds on $\omega' + i \faq^*$
for some compact
neighborhood $\omega'$ of $\omega$ and some polynomial
function $q' \in \Pi_\gS(\faqd),$
then it will hold on $\omega + i \faq^*$
for {\rm every} polynomial function $q \in \Pi_\gS(\faqd)$
for which $q(\gl) \nC_{Q|P}(s\col \gl)$ is holomorphic on a neighborhood
of $\omega + i \faqd$
(with constants $N$ and $C$ depending on $q$). This can be seen
by an application of the Cauchy integral formula as in \cite{Bps2},
proof of Lemma 6.1.
\end{rem}

\begin{rem}
Notice also that by \cite{BSft}, eq. (72),
the inverse of the normalized c-function is given by
\begin{equation}\naam{e: inverse C}
\nC_{Q|P}(s \col \gl)^{-1}=\nC_{P|Q}(s^{-1} \col s \gl).
\end{equation}
Hence
it follows from the estimate (\ref{e: Eqfa})
that this inverse satisfies
a similar bound on $\omega+i\faqd$,
although of course with possibly
different $q, N, C$,
\begin{equation}
\naam{e: Eqfa inverse}
\| q(\gl)\nC_{Q|P}(s\col\gl)^{-1} \| \leq C\, (1+|\gl|)^{N}.
\end{equation}
\end{rem}

Before we give the proof of
Thm. \ref{t: main thm}, we will review the relation
of the c-functions to standard intertwining
operators, meanwhile fixing useful notation.
We recall from the previous section
that every parabolic subgroup $P \in \minparabs$ has a Langlands decomposition
of the form (\refer{e: Langlands deco P min}), and that $\Mfu$ denotes the set of equivalence
classes of finite dimensional irreducible unitary representations of $M.$

If $\xi \in \Mfu,$ let $\cH_\xi$ be a Hilbert space in which $\xi$
is unitarily realized.
It is clear that $\Mfu\subset\cR(M)$.
We identify $\faqdc$ with the subspace of $\fadc$
consisting of elements vanishing on $\fa \cap \fh.$
We use the notation from Section \ref{s: estimates operators}
for the principal series representations $\pi_{P,\xi,\gl}$,
where $\xi\in\Mfu$ and $\gl\in\faqdc$. Note that
$\rho_P = \frac 12 \tr (\ad(\dotvar)|_{\fn_P})$ is defined
as a linear functional on $\fa,$ and that it
vanishes on $\fa \cap \fh,$ hence belongs to $\faqd;$
see \bib{Bps1}, Lemma 3.1.

Let $dk$ be normalized Haar measure on $K.$
The sesquilinear pairing
$$\Ci (P\col \xi \col \gl) \times \Ci (P \col \xi \col - \bar \gl) \to \C$$
given by
\begin{equation}
\naam{e: defi pairing}
\inp{f}{g} = \int_K \inp{f(k)}{g(k)}_\xi\; dk
\end{equation}
is $G$-equivariant. In the compact picture this
pairing becomes a pre-Hilbert structure on $\Ci(K \col \xi)$ for which
$\pi_{P, \xi, - \bar \gl}(x)^* = \pi_{P, \xi, \gl}(x)^{-1},$
for all $x \in G.$ In particular, the representation $\pi_{P, \xi, \gl}$
is unitarizable for $\gl \in i\faqd.$

Let $C^{-\infty}(K \col \xi)$ denote the space of
generalized functions $K \to \cH_\xi$, transforming according to the rule
(\refer{e: tr rule for K}), and for each $\gl\in\faqdc$ let
$C^{-\infty}(P \col \xi \col \gl)$ denote the space of generalized
functions $G \to \cH_\xi$ transforming according to the rule
(\refer{e: tr rule pr sr}).
These spaces are equipped with the usual locally convex
strong dual topologies of distribution spaces.
By (\refer{e: defi pairing})
we can identify $C^{-\infty}(K \col \xi)$
with the strong anti-linear dual of the Fr\'echet space
$C^\infty(K \col \xi),$
and likewise we can identify $C^{-\infty}(P \col \xi \col \gl)$
with the strong anti-linear dual of $C^{\infty}(P \col \xi \col -\bar\gl)$.
It follows by transposition that the restriction map $f \mapsto f|_K$
extends to a topological linear isomorphism
from $C^{-\infty}(P \col \xi \col \gl)$
onto $C^{-\infty}(K\col\xi)$.
The representation $\pi_{P, \xi, \gl}$ naturally extends to
a continuous representation $\pi_{P, \xi, \gl}$ of $G$
in both of these two model spaces. The pairing
(\ref{e: defi pairing}) is $G$-equivariant
with respect to these actions.

For $0<s<\infty$ let $C^{-s}(K \col \xi)$ denote the subspace of
$C^{-\infty}(K\col\xi)$ consisting of the generalized functions
of order at most $s$. The pairing (\refer{e: defi pairing})
extends to a  $K$-equivariant non-degenerate sesquilinear pairing
\begin{equation}
\label{pairing for distributions}
C^{-s}(K \col \xi) \times C^s(K \col \xi) \to \C.
\end{equation}
Recall that we have equipped $C^s(K \col \xi)$
with the norm $\|\dotvar\|_s$.
We equip $C^{-s}(K \col \xi)$
with the dual norm, denoted $\|\dotvar \|_{-s}.$ Then
\begin{equation}
\label{C-S inequality for distributions}
|\inp{f}{g}|\leq \|f\|_{-s}\|g\|_s
\end{equation}
for $f\in C^{-s}(K \col \xi)$ and $g\in
C^{s}(K \col \xi)$.
In this fashion
$C^{-s}(K \col \xi)$ becomes identified
with the anti-linear dual of the Banach space
$C^s(K \col \xi),$ as a representation space for $K$.

Let $P,Q \in \minparabs.$ For $R \in \R$ we define the set
$$
\faqd(Q\sepbar P, R): = \{\gl \in \faqdc \mid \inp{\Re \gl}{\ga} > R, \; \forall \ga \in \gS(\bar Q) \cap \gS(P)\},
$$
in analogy with (\ref{e: defi faPd Q sepbar P}).
\begin{lemma}
\naam{l: preparation intertwiner}
\begin{enumerate}
\itema
$\gS(P) = \{\ga |_{\faq}\; \mid \, \ga \in \gS(P, \fa)\,\},$
\itemb
$\Pi_{\gS(\bar Q) \cap \gS(P)}(\faqd)
= \{ p|_{\faqdc} \; \mid p \in \Pi_{\gS(\bar Q, \fa)\cap \gS(P, \fa)}(\fad)\;\},$
\itemc
for every $R \in \R,$ $\; \faqd(Q\sepbar P, R)= \fad(Q \sepbar P, R) \cap \faqd.$
\end{enumerate}
\end{lemma}
\proof
The nilpotent radical $\fn_P$ decomposes as the direct sum of the weight spaces for $\fa = \fa_P,$
with $\gS(P, \fa)$ as the associated set of weights. The set of the right-hand side of the equation in (a)
equals the set of $\faq$-weights in $\fn_P;$ this is $\gS(P).$ Assertion (a) follows.

We agreed to identify
$\faqd$ with the space of linear functionals in $\fad$ vanishing on $\faP \cap \fh.$
By symmetry of $\gs,$ the decomposition (\ref{e: deco fa}) is orthogonal. This implies that
the dual inner product on $\faqd$ coincides with the restriction of the dual inner product on
$\fad.$ Moreover, the restriction map $\eta \mapsto \eta|_{\faq}$ coincides with the
orthogonal projection $\fad \to \faqd.$
Assertions (b) and (c) now follow from (a).
\qed

It follows from Lemma \refer{l: preparation intertwiner} (c) that for $\gl \in \faqd(Q\mid P, r_\xi)$
the intertwining operator $A(Q \col P \col \xi \col \gl)$ is well defined by the convergent integral
(\refer{e: defining integral for A}). In view of Lemma \refer{l: preparation intertwiner} (b), (c),
the following result is now an immediate consequence of Theorem \refer{t: estimate A on smooth func}.

\begin{prop}
\naam{p: estimate A on smooth func, two}
Let $P, Q \in \minparabs.$

For every $f \in \Ci(K\col \xi)$ the function
$\gl \mapsto A(Q\col P \col \xi \col \gl) f,$
originally defined on $\faqd(Q|P, r_\xi),$
extends to a meromorphic $\Ci(K\col \xi)$-valued
function on $\faqdc.$

Let $R \in \R.$
There exist  a polynomial function
$q\in\Pi_{\gS(\bar Q) \cap \gS(P)}(\faqd)$ and constants
$r\in\N$ and $N \in \N,$
such that the following is valid.

The map
$\gl \mapsto q(\gl)A(Q\col P\col\xi\col\gl)f$
is holomorphic as a $C^\infty(K\col \xi)$-valued function on
$\faqd(Q |P,R)$, for
every $f \in \Ci(K\col \xi)$. Moreover, for
each $s \in \N$ there exists a constant $C > 0$ such that
\begin{equation}
\naam{e: estimate A on smooth func}
\|q(\gl)A(Q \col P\col\xi\col\gl) f\|_{s}\leq C(1+|\gl|)^N\, \|f\|_{s + r},
\end{equation}
for all $\gl\in\faqd(Q | P,R)$ and $f \in \Ci(K\col \xi)$.
\end{prop}

\begin{rem}
\naam{r: estimate A, two}
It follows by the above estimate that, for
$\gl \in \faqd(Q|P, R),$
the map $A_0(\gl):= q(\gl) A(Q \col P \col \xi \col \gl):$
$C^\infty(K\col \xi) \to C(K\col \xi)$
extends to a bounded linear map from
$C^{s + r}(K \col \xi)$ to $C^s(K\col \xi),$ for every $s \in \N.$

Moreover, by a standard application of the Cauchy estimates
for power series coefficients,
we infer, with $R, $ $r$ and $q$ as above,
that $\gl \mapsto A_0(\gl)$
is holomorphic on $\faqd(Q|P, R)$ as a function with values in the Banach space
$B_{s, s+ r}$ of bounded linear maps
$C^{s+r}(K\col \xi) \to  C^{s}(K \col \xi).$
\end{rem}

With respect to the pairing (\refer{e: defi pairing}) it is known that
\begin{equation}
\naam{e: adjoint of A}
A(Q \col P \col \xi \col \gl) = A(P \col Q\col \xi \col -\bar \gl)^*,
\end{equation}
for generic $\gl \in \faqdc.$ In fact, this follows from the validity of this
identity
on the space $C^\infty(K\col \xi)_K$ of $K$-finite functions in
$C^\infty(K\col \xi);$ see
\cite{KS2}, Prop.\ 7.1 (iv).

It follows that the intertwining operator (\refer{e: standard intertwiner})
admits a continuous linear extension to an intertwining operator
$$
A(Q \col P \col \xi \col \gl): \;\; \Cminf(P\col \xi \col \gl) \;\too \;
\Cminf(Q\col \xi \col \gl),
$$
given by the formula (\refer{e: adjoint of A}), for generic $\gl \in  \faqdc.$

The space of generalized functions $f\in \Cminf(P \col \xi \col \gl)$
which are invariant under $\pi_{P,\xi\gl}(H)$ is denoted by
$\Cminf(P \col \xi \col \gl)^H$. We need the following
description of this space from \bib{Bps1}.

Let $\xi\in\Mps$. For $v \in \cW$ we define $V(\xi, v)$ to be the space
$\cH_\xi^{M \cap v H v^{-1}},$
equipped with the restricted Hilbert inner product. Moreover, we define
the space $V(\xi)$ to be the formal orthogonal direct sum
$$
V(\xi) := \oplus_{v \in \cW}\;\;  V(\xi, v).
$$
Note that $V(\xi)$ is non-trivial by the definition of $\Mps$,
see Section \ref{s: eisenstein integrals}.

{}From \bib{Bps1}, Sect.\ 5, we recall the
definition of the evaluation
map $$\ev : \Cminf(P \col \xi \col \gl)^H \to V(\xi)$$ by
$$
[\ev (\gf)]_v = \gf(v), \qquad (v \in \cW).
$$
This map is well defined, as an $H$-fixed element of $\Cminf(P \col \xi \col \gl)$
is smooth on each open orbit $P v H,$ for $v \in \cW.$
According to \bib{Bps1}, Thm.\ 5.10, there exists
a unique meromorphic $\Hom(V(\xi), C^{-\infty}(K\col \xi))$-valued function
$j(P\col \xi \col \dotvar)$ on $\faqdc,$ such that for every
$\eta \in V(\xi)$ and generic $\gl\in\faqdc$,
\begin{enumerate}
\vspace{4pt}
\itema
the generalized function $j(P \col \xi \col \gl)\eta$
belongs to $\Cminf(P \col \xi \col \gl)^H$;
\itemb
$\ev \after j(P \col \xi \col \gl)\eta = \eta.$
\end{enumerate}

We also need the operator $B(Q \col P \col \xi \col \gl)\in\End(V(\xi))$
constructed in \bib{Bps1}, Sect.\ 6.
Recall that this depends meromorphically
on $\gl\in\faqdc$, and that it is
uniquely determined by
the relation
\begin{equation}
\naam{e: formula for B}
B(Q \col P \col \xi \col \gl)  =
\ev \after A(Q \col P \col \xi \col \gl) \after j(P \col \xi \col \gl),
\end{equation}
for generic $\gl \in \faqdc.$

The space $\oC = \oCtau$ admits a  direct sum decomposition of the form
\begin{equation}
\naam{e: deco oCtau}
\oCtau = \oplus_{\xi \in \Mps} \;\; \oC_\xi(\tau),
\end{equation}
as in \bib{BSft}, proof of Lemma 3.
All but finitely many summands of this direct sum are trivial.

Let
$$
C(K \col \xi \col \tau): =  [C(K \col \xi) \otimes \Vtau]^K.
$$
Then according to \bib{BSft}, Lemma 3,
there exists, for each $\xi \in \Mps,$
a natural linear isomorphism
$T \mapsto \psi_T$ from $C(K\col\xi\col\tau) \otimes \overline{V(\xi)}$
onto $\oCtau_\xi,$ where $\overline{V(\xi)}$
is the conjugate Hilbert space of $V(\xi)$.

Via these isomorphisms,  the $c$-function
$\nC_{Q|P}(1\col\gl)$
may be expressed in terms of standard intertwining operators, as follows.
For
$ f \otimes \eta \in C(K\col\xi\col\tau) \otimes \overline{V(\xi)},$
we have
\begin{equation}
\naam{e: C in terms of A B}
\nC_{Q|P}(1\col\gl)\psi_{f\otimes\eta}=
\psi_{A(Q\col P\col\xi\col -\gl)f
\otimes B(\bar P\col\bar Q\col\xi\col\bar\gl)^{-1}\eta},
\end{equation}
as an identity of meromorphic functions in the variable $\gl \in \faqdc;$
see \bib{BSft}, Eqn.\ (57).
In particular, it follows that
$\nC_{Q|P}(1 \col \gl)$ preserves the decomposition (\refer{e: deco oCtau}).
It is convenient to rewrite (\ref{e: C in terms of A B})
in the following form, which follows
by application of (\refer{e: inverse C})
\begin{equation}
\naam{e: C in terms of A B, two}
\nC_{Q|P}(1\col\gl)\psi_{f\otimes\eta}=
\psi_{A(P\col Q\col\xi\col -\gl)^{-1}f
\otimes B(\bar Q\col\bar P\col\xi\col\bar\gl)\eta}.
\end{equation}
Here we note that $\gl\mapsto B(\bar Q\col\bar P\col\xi\col\bar\gl)\eta$ is a meromorphic
$\overline{V(\xi)}$-valued function.

We now come to the proof of the theorem.
\medbreak
\noindent{\it Proof of Theorem \refer{t: main thm}. \ }
It follows from the relation
\bib{BSft}, Eqn.\ (68), that we may assume that $s=1$.
Thus, we aim to prove:
\begin{equation}
\naam{e: Eqfb}
\|q(\gl)\nC_{Q|P}(1 \col\gl)\|\leq C(1+|\gl|)^{N},\qquad
(\gl\in\omega+i\faqd).
\end{equation}

Furthermore, it follows from (\ref{e: C in terms of A B})
combined with the argument
below Lemma 3.2 in \cite{BSfi}, that for every compact
subset $\omega' \subset \faq^*$
there exists a function $q' \in \Pi_{\gS, \R}(\faqd)$ such that
$q'(\gl) \nC_{Q|P}(1 \col \gl)$
is holomorphic on a neighborhood of $\omega' + i \faq^*.$
By combining this observation with that of
Remark \ref{r: some implies all}, we see that
it suffices to establish (\ref{e: Eqfb}) with some $q \in
\Pi_\gS(\faqd).$

{}From (\refer{e: deco oCtau}) and (\refer{e: C in terms of A B}),
and the text between and after these displays,
it follows that it suffices
to establish, for every $\xi \in \Mps,$
the existence of $q,$ $N$ and $C$ such that
the estimate (\refer{e: Eqfb}) holds with $\nC_{Q|P}(1 \col \gl)$ replaced
by its restriction $ \nC_{Q|P}(1 \col \gl)_\xi$ to $\oC_\xi(\tau).$ We fix
$\xi \in \Mps.$

It follows from \bib{Bps2}, Lemma 16.6,
that for every $R \in \R$ there exist a polynomial $q\in\Pi_{\gS}(\faqd)$ and
constants $N$, $C$ such that
$$
\|q(\gl)
A(P\col Q\col\xi\col -\gl)^{-1}f\|_0\leq C(1+|\gl|)^N\|f\|_0
$$
for all $\gl\in\faqd(Q\sepbar P,R)$ and $f\in C(K\col\xi\col\tau)$.
In particular, this can be achieved for
$\gl\in\omega+i\faqd$ by an appropriate
choice of $R$.
It thus only remains to obtain a similar bound for
$B(\bar Q\col\bar P\col\xi\col\bar\gl)\eta$. This is accomplished in
Proposition \ref{p: polynomial bound B}
of the next section.
\qed

\eqsection{Estimates for $B(Q:P:\xi:\gl)$}
\label{s: estimates B}
In the following it will be convenient to use the following
notation for certain subsets of $\faqdc,$ for given
$P,Q \in \minparabs$ and $R \in \R:$
$$
\cA(Q,P,R) :=
\{\gl\in \faqdc  \mid  \forall \ga\in\gS(\bar Q)\cap\gS(P):\;
|\inp{\Re\gl}{\ga}|< R\}.
$$
We also define
$$\faqd(P,R):=\{\gl \in \faqdc \mid
\inp{\Re \gl}{\ga} < R, \; \forall \ga \in\gS(P)\},
$$
and observe that in particular
$$\cA(\bar P,P,R)=\faqd(P,R)\cap\faqd(\bar P,R).$$
Finally, we note that for each compact set $\omega\subset\faqd$
there exists $R>0$ such that $\omega+i\faqd\subset
\cA(Q,P,R)$.

\begin{prop}\naam{p: polynomial bound B}
Let $P, Q\in\Parmin$, $\xi\in\Mps$ and $R>0$. There exists
a polynomial function $q \in \Pi_{\gS}(\faqd)$ and constants $N \in \N$ and $C > 0$
such that the following is valid. The map
$
\gl \mapsto q(\gl)B(Q\col P\col\xi\col\gl)
$
is a holomorphic
$\End(V(\xi))$-valued function on the set $\cA(Q,P, R).$
Moreover, for all $\gl$ in this set,
$$
\|q(\gl)B(Q\col P\col\xi\col\gl)\| \leq C(1+|\gl|)^N.
$$
\end{prop}

\proof
By the argument in \bib{BSmc}, Lemma 20.5, we may assume
that $Q= \bar P$.
The condition on $\gl$
is then that $\gl\in\faqd(P,R)\cap\faqd(\bar P,R)$.

For generic $\gl$ the endomorphism $B(\bar P \col P\col\xi\col\gl)$
of $V(\xi)$ is given as the composition of three maps
in (\refer{e: formula for B}). In terms of the compact
picture we will view these as maps
\begin{eqnarray}
&&j(P\col \xi\col \gl):\; \; V(\xi) \;\to \; \Cminf(K\col \xi),\naam{e: map j}\\
&&A(\bar P \col P \col \xi \col \gl):\;\;
       \Cminf(K \col \xi)\;\to\; \Cminf(K \col \xi),\naam{e: map A}\\
&&\ev:\; \; \Cminf(K \col \xi)^{H}_{\bar P, \gl} \;\to \; V(\xi).
\naam{e: map ev}
\end{eqnarray}
The indices $\bar P , \gl$ on the domain of $\ev$ indicate that the space
of $H$-invariants for the representation $\pi_{\bar P ,\xi, \gl}$ has
been taken.

It will be seen in the following
Lemmas \refer{l: polynomial bound j}, \refer{l: polynomial bound A on Cminf} and \refer{l: polynomial bound for ev}
that each of the three maps above
satisfies a bound of the desired polynomial type in $\gl,$ in terms
of suitable operator norms. The estimates will turn out to be valid on sets of the
form $\faqd(P,R_1),$
$\faqd(\bar P| P, R_2)$ and $\faqd(P,R_3),$ respectively, with $R_j$
arbitrary real numbers. Noting that
$\faqd(\bar P| P, R_2) = \faqd(\bar P, -R_2)$
and taking $R_1 = R_3 = R$ and $R_2 = -R$ we obtain sets whose intersection equals
$\faqd(P, R) \cap \faqd(\bar P, R).$ Thus, combining the mentioned lemmas,
the proof of Proposition \refer{p: polynomial bound B} is completed.
\qed

\begin{lemma}
\naam{l: polynomial bound j}
Let $P\in\Parmin$, $\xi\in\Mps$ and
$R\in\R$. There exist $s\in\N^+$ and $q\in\Pi_{\gS}(\faqd)$
such that the restriction to $K$ of $q(\gl)j(P\col\xi\col\gl)\eta$
belongs to $C^{-s}(K\col\xi)$
for all $\eta\in V(\xi)$ and $\gl\in\faqd(P,R)$,
and such that the norm of
the restriction satisfies the following uniform estimate. There exist
constants $N\in\N$, $C>0$ such that
$$\|q(\gl)j(P\col\xi\col\gl)\eta\|_{-s}\le C(1+|\gl|)^N\|\eta\|$$
for all $\eta\in V(\xi)$, $\gl\in\faqd(P,R)$.
\end{lemma}

\proof
See \bib{Bps2}, Thm.\ 9.1.
\qed

We now turn to estimating a suitable operator norm
of the intertwining operator (\refer{e: map A}). As this map
is given as an adjoint in (\refer{e: adjoint of A}), the following
lemma is the key step.
In what follows, we shall use, for $p,q \in \Z,$
the following abbreviation for the space
of bounded linear maps between Banach spaces:
$$
B_{q,p}:= B(C^p(K\col \xi), C^q(K\col \xi)).
$$
Moreover, we shall use the notation $\|\dotvar\|_{q,p}$ for the operator
norm on this space.
Accordingly, if $T \in B_{q,p}$ and $T' \in B_{r,q},$ then $T'\after T \in B_{r,p}$
and
$$
\|T' \after T \|_{r,p} \leq \|T'\|_{r,q} \|T\|_{q,p}.
$$

\begin{lemma}
\naam{l: polynomial bound A on Cminf}
Let  $P\in\Parmin$, $\xi\in\Mfu$ and
$R\in\R$.  There exist a polynomial function $q\in\Pi_{\gS}(\faqd),$
and constants  $r\in \N$ and $N \in \N$ such that the following holds.

For every $s \in \N\setminus\{0\}$ there exists a constant $C > 0$
such that for every
$\replacefbythis \in C^{-s}(K\col \xi),$
the map $\gl \mapsto q(\gl)A(\bar P\col P\col\xi\col\gl)\replacefbythis$
is holomorphic as a $C^{-r-s}(K\col\xi)$-valued function
on $\faqd(\bar P | P, R)$
and
\begin{equation}
\naam{e: estimate A on gen func}
\|q(\gl)A(\bar P\col P\col\xi\col\gl) \replacefbythis \|_{-r-s}
\leq C(1+|\gl|)^N\,\ \| \replacefbythis\|_{-s} ,
\end{equation}
for all  $\gl\in\faqd(\bar P | P, R)$.
\end{lemma}

\proof
We will derive the result from Proposition
\refer{p: estimate A on smooth func, two} by using
(\refer{e: adjoint of A}).

For $q \in \Pi_\gS(\faqd),$ we put $q^\vee(\gl): = \overline{q(- \bar \gl)}.$
Then $q \mapsto q^\vee$ is a bijection of $\Pi_{\gS}(\faqd)$ onto itself.

Let $s, t$ be positive integers,
and let $T \in B_{t,s}.$ Then we define the conjugate map
$T^* : C^{-t}(K\col \xi) \to C^{-s}(K \col \xi)$ by $\inp{T^* \gf}{\replacefbythis} =
\inp{\gf}{T \replacefbythis},$ for all $\gf \in C^{-t}(K\col \xi)$ and
$\replacefbythis \in C^s(K\col \xi).$ It is readily seen that $T \mapsto T^*$
defines an anti-linear map $B_{t,s} \to B_{-s, -t},$
with
$$
\|T^*\|_{-s, -t} \leq \|T\|_{t,s}
$$
for all $T \in B_{t,s}.$

Let now $R \in \R.$ Let  $q_0 \in \Pi_\gS(\faqd),$ $r \in \N$
and $N \in \N$ be as in Proposition \refer{p: estimate A on smooth func, two}
with $(P, \bar P)$ in place of $(Q, P).$
Let $s \in \N$ and let $C > 0$ be a constant as in
(\refer{e: estimate A on smooth func}).
According to Remark \refer{r: estimate A, two}
it follows from the mentioned proposition
that the map $T: \gl \mapsto q_0(\gl) A(P \col \bar P \col \xi \col \gl)$
is a holomorphic function on $\faqd(P|\bar P, R)$
with values in $B_{s, s + r},$
satisfying $\|T(\gl)\|_{s, s+r} \leq C (1 + |\gl|)^N$ for all
$\gl \in \faqd(P|\bar P, R)$.
This implies that, for every $s \in \N,$ the
map $\gl \mapsto T(- \bar \gl)^*$ is a holomorphic
function with values in $B_{-s - r, -s},$ satisfying
the estimate
$$
\|T(-\bar \gl)^*\|_{-s - r, -s} \leq \|T(- \bar \gl)\|_{s, s + r}
\leq C (1 + |\gl|)^N,
$$
for all $\gl \in - \faqd(P|\bar P, R).$ We now observe that
$- \faqd(P|\bar P, R) = \faqd(\bar P|P, R),$ and that
$T(-\bar \gl)^* = q^\vee_0(\gl) A(\bar P \col P \col \xi \col \gl)$
in view of (\refer{e: adjoint of A}).
The estimate (\refer{e: estimate A on gen func}) thus follows with
$q = q_0^\vee.$
\qed

\begin{lemma}
\naam{l: polynomial bound for ev}
Let
$P\in\Parmin$, $\xi\in\Mps$, $R\in\R$ and $s\in\N$. There exists a constant $C>0$ such that
\begin{equation}
\label{e: estimate ev f}
\|\ev \;\replacefbythis\|\le C(1+|\gl|)^s \|\replacefbythis\|_{-s}
\end{equation}
for all $\gl \in \faqd(\bar P, R)$ and
$\replacefbythis\in
C^{-s}(K\col\xi)\cap C^{-\infty}(K\col\xi)^H_{P,\gl}.$
\end{lemma}

The proof of this lemma will be given in the next section.

\eqsection{Estimation of the evaluation map}\label{s: Esteval}
This section is entirely devoted to the proof of Lemma
\ref{l: polynomial bound for ev}.
To prepare for it, we first make the following general observation.

Let $G$ be a Lie group, and let $P$ be a closed subgroup.
Then the adjoint action of $P$ on $\fg$ naturally
induces a left action of $P$ on the space of densities on
$\fg/\fp.$  As the latter space is one-dimensional, the action
is given by a character $\chi: P\to \R^+.$
In fact,
 this character is given by
$$
\chi(p) = \frac{|\det \Ad_P(p)|}{|\det \Ad_G(p)|}, \qquad (p \in P).
$$
We consider the space $C(G\col P\col \chi)$ of continuous functions
$f: G \to \C$ such that $f(p g) = \chi(p) f(g)$ for all
$g \in G$ and $p \in P$ and write $C_c(G\col P\col \chi)$ for the
subspace of functions with compact support modulo $P.$

Let $\cD_{P \bs G}$ denote the density bundle on $P \bs G.$ Its fiber at a point $Pg$ is the space $\cD(T_{Pg}(P\bs G))$ of densities on the tangent space $T_{Pg} (P\bs G).$ Pull-back by right multiplication gives
a linear isomorphism $d(r_{g^{-1}})(Pg)^*: \cD(T_{Pe}(P\bs G)) \to
\cD(T_{Pg}(P\bs G)).$

Let $\Gamma_c( \cD_{P \bs G})$ denote the  space of compactly supported continuous sections of $\cD(P\bs G).$
We fix a positive density $\omega$ on $T_{Pe}(P\bs G)\simeq \fg/\fp.$
Then the
map $f \mapsto j_\omega(f)$ given by
$$
j_\omega(f)(Pg) = f(g)\; d(r_{g^{-1}})(Pg)^* \omega
$$
defines a topological linear isomorphism
$$
j_\omega: C_c(G \col P \col \chi) \;\;{\buildrel \simeq \over \longrightarrow}\;\;
\Gamma_c(\cD_{P\bs G}).
$$

\begin{lemma}
\label{l: about density}
Let $L$ be a closed subgroup of $G$ such that $PL$ is open in $G$
and such that $\chi|_{L \cap P} = 1.$ Then there exists a unique right-invariant positive Radon measure $d\bar l$ on $(L\cap P) \bs L$ such that
$$
\int_{P\bs G} j_\omega(f) = \int_{(L\cap P )\bs L} f(l) \; d\bar{l}
$$
for all $f \in C_c(G\col P \col \chi)$ with $\supp f \subset PL.$
\end{lemma}

\proof
Let $C_L(G\col P\col \chi)$ denote the space of functions $f$ as in the lemma.
Then $f \mapsto f|_L$ defines a topological linear isomorphism
from $C_L(G\col P\col \chi)$ onto $C_c((L\cap P) \bs L).$ We denote its
inverse by $g \mapsto \tilde g.$
The functional $\mu: g \mapsto \int_{P\bs G} j_\omega(\tilde g)$ defines a Radon measure on $(P\cap L)\bs L.$ Let $m \in L.$ As the integration of densities is invariant under pull-back by diffeomorphisms, we find,
for $g \in C_c((L\cap P)\bs L),$ that
$$
\mu(r_m^* g) = \int_{P\bs G} j_\omega(\widetilde {r_m^* g}) =
\int_{P\bs G} r_m^* (j_\omega(\tilde {g})) = \int_{P\bs G} j_\omega(\tilde g)= \mu(g).
$$
Hence,  $\mu$ is right $L$-invariant and the result follows with
$d\bar l = \mu.$ Uniqueness is obvious.
\qed

We now return to the notation of the previous sections,
with $G$ a group of the Harish-Chandra class,
and apply Lemma \ref{l: about density} with
$P$ a parabolic subgroup
from $\Parmin.$  Then
$$
\chi(man) = a^{2\rho_P}, \qquad (m \in M_P, a \in A_P, n \in N_P).
$$
It follows that $\chi = 1$ on both $K \cap P = K \cap M_P$
and $H \cap P = (H\cap M_P)A_H.$ So we may apply Lemma \ref{l: about density} with $L$ equal to $K$ and with $L$ equal to $H.$

It follows that we may fix the density $\omega$ on $\fg/\fp$ (uniquely) such that the associated invariant
measure $d \bar k$ on $(M_P\cap K) \bs K$ is normalized.
Furthermore, let  $d \bar h$ be the invariant measure on $(H\cap P) \bs H$
determined by Lemma \ref{l: about density} with $L = H.$

\begin{cor}\label{c: Olafssons lemma}
For every $f\in C(G\col P\col 1\otimes 2\rho_P\otimes 1)$ with
$\supp f \subset PH,$ we have
$$
\int_{(K\cap M_P)\bs K} f(k) \; d\bar k = \int_{(H \cap P)\bs H } f(h) \; d \bar h.
$$
\end{cor}

\proof
By application of Lemma \ref{l: about density}, both integrals equal the integral
of the density $j_\omega(f)$ over $P\bs G.$
\qed

\begin{rem}
The above result is due to \bib{Olfp}, Lemma 1.3.
Our proof is more conceptual.
\end{rem}

\noindent%
{\em Proof of Lemma \ref{l: polynomial bound for ev}.\ }
It suffices to prove an estimate similar to
(\ref{e: estimate ev f}) for each of the finitely
many components $\ev_w \replacefbythis$ of $\ev \replacefbythis,$ for $w \in \cW.$ Let
$R_w$ denote the topological linear
automorphism of $C^{-\infty}(K\col \xi)$
given by $R_w\replacefbythis(k) = \replacefbythis(kw).$
Then $R_w$ maps $C^{-s}(K\col\xi)$ isomorphically to itself for each
$s>0$, and it
restricts to a bijection
$$
C^{-\infty}(K \col\xi)^H_{P, \gl}\;\;
{\buildrel \simeq \over \longrightarrow}
\;\;C^{-\infty}(K \col\xi)^{wHw^{-1}}_{P, \gl}.
$$
Since
$\ev_e \after R_w = \ev_w,$ it suffices to prove the estimate (\ref{e: estimate ev f})
with $\ev_e$ in place of $\ev.$

Fix an element $v\in\cH_\xi^{M\cap H}$.
Then we shall complete the proof by establishing the
existence of $C>0$ such that
\begin{equation}
\naam{e: Eqfe}
|\hinp{\replacefbythis(e)}v_\xi|\le C(1+|\gl|)^s \|\replacefbythis\|_{-s},
\end{equation}
where $\hinp{\,\cdot\,}{\,\cdot\,}_\xi$ denotes the unitary structure
on $\cH_\xi$.  Fix
a non-trivial and non-negative compactly supported smooth function $\psi$
 on $(H\cap P)\bs H.$
Define, for each $\gl\in\faqdc$, an
$\cH_\xi$-valued function $\gf_\gl$ on $G$ as follows:
\begin{equation}
\naam{e: Eqfi}
\gf_\gl(x)=
\cases{
a^{-\bar\gl+\rho}\psi(h)\xi(m)v & if $x=manh\in PH,$\cr
0 & otherwise,\cr}
\end{equation}
where $m\in M$, $a\in\Aq$, $n\in N_P$, $h\in H$. Then, since $\psi$ is
compactly supported, the support of $\gf_\gl$ is a closed subset of
the open subset $PH$ of $G$, and hence
$\gf_\gl\in C^\infty(P\col\xi\col-\bar\gl)$ (cf.\ \bib{HSbook}, Thm.\ II.3.3).
We claim that, with respect to the sesquilinear pairing
(\ref{pairing for distributions}),
\begin{equation}
\naam{e: Eqff}
\hinp \replacefbythis{\gf_\gl}=
\hinp{\replacefbythis(e)}v_\xi\int_{(H\cap P)\bs H} \psi(h) \, d\bar h
\end{equation}
for all $\replacefbythis\in C^{-\infty}(P\col\xi\col\gl)^H.$
Indeed, since $\replacefbythis$ is smooth
on $PH$, which contains the support of $\gf_\gl$, we have
$$
\hinp \replacefbythis{\gf_\gl}= \int_{K} \hinp{\replacefbythis(k)}{\gf_\gl(k)}_\xi \,dk,
$$
and by Corollary \ref{c: Olafssons lemma}, the integral over $K$ can be rewritten
as the following integral over $(H\cap P)\bs H$:
$$
\int_{(H\cap P)\bs H} \hinp{\replacefbythis(h)}{\gf_\gl(h)}_\xi \,d\bar h.
$$
Since $\replacefbythis(h)=\replacefbythis(e)$ and $\gf_\gl(h)=\psi(h)v$,
the claimed
formula (\refer{e: Eqff}) follows.

The integral in (\refer{e: Eqff}) is positive. Applying
the inequality (\ref{C-S inequality for distributions}) to
the left side of
(\ref{e: Eqff}) we infer that
$$
|\hinp{\replacefbythis(e)}v_\xi|\le C\,\|\replacefbythis\|_{-s}\|\varphi_\gl\|_s
$$
for some constant $C$.
We will finish the proof
of (\ref{e: Eqfe}) by showing that given $R \in \R$ and $s\in\N$,
there exists $C>0$ such that
\begin{equation}
\naam{e: Eqfg}
\|\gf_\gl\|_s\le C(1+|\gl|)^s
\end{equation}
for all $\gl\in\faqd(\bar P,R)$.

Fix a compact set $\Omega\subset H$ whose image contains $\supp\psi$
under the mapping $h\mapsto (H\cap P)h$, $H\to (H\cap P)\bs H$.
Then $\supp\gf_\gl\subset P\Omega$ for all $\gl$.
If $\fb$ is a Lie algebra, and $s \in \N,$ then by
$U(\fb)_s$ we denote the subspace of elements in $U(\fb)$ of order
at most $s$.
The norm $\|\gf_\gl\|_s$ is dominated by a positive constant times the maximum
of $\sup_{k\in K}\| R_u\gf_\gl(k)\|_\xi$
for $u$ ranging over a finite subset $S$ of $U(\fk)_s$
(see (\ref{d: s-norm})).
Taking into account
that $- \faqd(P, R) = \faqd(\bar P, R),$ we see from Lemma
\refer{l: sup on K and H} below that we may as well estimate the norm
$\sup_{h\in\Omega}\|(\pi_{P, \xi, - \bar \gl}(u)\gf_\gl)(h)\|_\xi$
for each  $u \in S.$
The estimate (\refer{e: Eqfg}) now follows by application of  Lemma \refer{l: estimate on Omega}, also given below.
\qed

\begin{lemma}
\naam{l: sup on K and H}
Let $P \in \Parmin$ and $\xi \in \Mfu.$
Let $\Omega\subset H$ be a compact set, and $R\in\R$.
There exists a constant $C>0$ such that
$$
\sup_{k\in K}\|\phi(k)\|_\xi\le C\sup_{h\in \Omega}\|\phi(h)\|_\xi,
$$
for all $\gl\in\faqd(P,R)$ and all
$\phi\in C(P\col\xi\col\gl)$
 with $\supp\phi\subset P\Omega.$
\end{lemma}

\proof We first give the proof under the assumption that $H$ is connected.
For $x\in PH$ we write $x=n(x)a(x)m(x)h(x)$ with
$n(x)\in N_P$, $a(x)\in\Aq$, $m(x)\in M$ and $h(x)\in H$.
The element $a(x)$ is unique and depends continuously on $x$.
Moreover, we may arrange that $h(x)\in\Omega$ for $x\in P\Omega.$
Then $\phi(k)=a(k)^{\gl+\rho}\xi(m(k))\phi(h(k))$ for $k\in K\cap PH$.
As $\xi$ is a unitary representation of $M$, it suffices
to prove that $a(k)^{\gl+\rho}$ is uniformly bounded for
$k\in K\cap P\Omega$ and $\gl\in\faqd(P,R)$.

For each $\ga\in\gS$, let $H_\ga\in\faq$ be determined by
$\nu(H_\ga)=\inp{\nu}{\ga}$ for all $\nu\in\faqd$.
It follows from \bib{Bconv}, Thm.\ 3.8,
that, for all $k \in K \cap PH,$ the element $\log a(k)$
belongs to the closed cone spanned by the
vectors $H_\ga$, for $\ga\in\gS(P).$
Moreover, by continuity and compactness, $\log a(k)$ belongs to
a bounded subset of this cone, for all $k\in K\cap P\Omega$.
Hence, there exist $c>0$ and numbers $r_\ga(k)\in [0,c]$ such that
$\log a(k)=\sum_{\ga\in\gS(P)} r_\ga(k) H_\ga$ for all $k\in K\cap P\Omega$.
It follows that
$$
(\Re\gl+\rho)(\log a(k)) = \sum_{\ga\in\gS(P)}
r_\ga(k) \inp{\Re\gl+\rho}{\ga}
$$
is bounded from above, uniformly for
$k\in K\cap P\Omega$ and $\gl\in\faqd(P,R)$.

It remains to treat the case where $H$ is not connected.
Let $H_e$ denote the identity component of $H.$ Then
$H=H_e \,(K \cap H).$
Hence, $\Omega$ can be written as a disjoint union
$\Omega_1k_1\cup\dots\cup \Omega_mk_m$, where  $k_i \in K \cap H,$
and where $\Omega_i\subset H_e$ is compact, for
$i=1,\dots,m$. Without loss of generality we may assume that
$\supp\phi\subset P\Omega_ik_i$ for some $i$, and then the result
follows by application of the above to the right translate of $\phi$ by $k_i^{-1}$.
\qed

\begin{lemma}
\naam{l: estimate on Omega}
Let $P \in \Parmin$ and $\xi \in \Mfu.$
Let $\gf_\gl \in C^\infty(P\col \xi \col \gl)$ be a
family of functions, for $\gl \in \faqdc,$ such that $\gf_\gl|_H$ is independent
of $\gl$ and has support inside a set of the form $(H \cap M) \Omega,$ with  $\Omega \subset H$ compact.
Then for every $s \in \N$ and every $u \in U(\fg)$ of order at most $s$ there exists
a constant $C >0$ such that
\begin{equation}
\label{e: estimate pi gf gl}
\sup_{h \in \Omega} \|\pi_{P, \xi, \gl} (u) \gf_\gl(h)\|_\xi \leq C\, (1 + |\gl|)^s,
\end{equation}
for all $\gl \in \faqdc.$
\end{lemma}

Let $h\in\Omega$. Since $\fg=\Lie(P)+\fh$ it follows from the
Poincar\'e--Birkhoff--Witt theorem that $u \in U(\fg)_s$ can be written as a finite sum of
products of the form $u_h:= \Ad h^{-1}(u'_h)u''_h$, where $u'_h\in U(\Lie(P))_s$
and $u''_h\in U(\fh)_s$. Moreover, since $\Omega$ is compact, the
elements $u'_h, u''_h$ can be chosen such that they
belong to bounded subsets of $U(\Lie(P))_s$ and $U(\fh)_s$,
respectively, for all $h\in\Omega$.

For each element $u'\in U(\Lie(P))_s$ there exists
a constant $C'>0$ such that
$$
\|(\xi \otimes (\gl + \rho_P) \otimes 1)(u')\|_{{\rm op}} \leq C'\, ( 1 + |\gl|)^s,
$$
where $\|\dotvar\|_{\rm op}$ indicates the operator norm on $\End(\cH_\xi).$
Moreover, the constant $C'$ can be taken  uniform
if $u'$ varies in a bounded subset of $U({\rm Lie(P)})_s.$

We apply this with $u'=u'_h.$ Then
$$
\|\pi_{P, \xi, \gl} (u_h) \gf_\gl(h)\|_\xi \leq C' (1 + |\gl|)^s\, \|\pi_{P, \xi, \gl}(u''_h)\gf_\gl(h)\|_\xi.
$$
It follows from the hypotheses that
$\|(\pi_{P,\xi,\gl}(u''_h)\gf_\gl)(h)\|_\xi$ is independent
of $\gl$ and uniformly bounded with respect to $h\in\Omega$.
The estimate (\ref{e: estimate pi gf gl})
now follows.
\qed

\def\adritem#1{\hbox{\small #1}}
\def\distance{\hbox{\hspace{3.5cm}}}
\def\apetail{@}
\def\adderik{\vbox{
\adritem{E.~P.~van den Ban}
\adritem{Mathematical Insitute}
\adritem{Utrecht University}
\adritem{PO Box 80 010}
\adritem{3508 TA Utrecht}
\adritem{The Netherlands}
\adritem{E-mail: E.P.vandenBan{\apetail}uu.nl}
}
}
\def\addhenrik{\vbox{
\adritem{H.~Schlichtkrull}
\adritem{Department of Mathematical Sciences}
\adritem{University of Copenhagen}
\adritem{Universitetsparken 5}
\adritem{2100 K\o benhavn \O}
\adritem{Denmark}
\adritem{E-mail: schlicht@math.ku.dk}
}
}
\vfill
\hbox{\vbox{\adderik}\vbox{\distance}\vbox{\addhenrik}}

\end{document}